
\input amstex.tex
\documentstyle{amsppt}
\magnification=\magstep1
\hsize=12.5cm
\vsize=18cm
\hoffset=1cm
\voffset=2cm
\def\DJ{\leavevmode\setbox0=\hbox{D}\kern0pt\rlap
{\kern.04em\raise.188\ht0\hbox{-}}D}
\footline={\hss{\vbox to 2cm{\vfil\hbox{\rm\folio}}}\hss}
\nopagenumbers 
\font\ff=cmr8
\def\txt#1{{\textstyle{#1}}}
\baselineskip=13pt
\def\hf{{\textstyle{1\over2}}}
\def\a{\alpha}
\def\d{{\,\roman d}}
\def\e{\varepsilon}
\def\f{\varphi}
 \def\G{\Gamma}

\def\s{\sigma}

\def\={\;=\;}
\def\zx{\zeta(\hf+ix)}
\def\zt{\zeta(\hf+it)}

\def\D{\Delta}  
\def\no{\noindent}  
\def\R{\Re{\roman e}\,}  \def\s{\sigma}
\def\z{\zeta} 
\def\no{\noindent} 
\def\e{\varepsilon}
\def\D{\Delta} 
\def\no{\noindent} 
\def\e{\varepsilon}

\def\no{\noindent} 
\font\teneufm=eufm10
\font\seveneufm=eufm7
\font\fiveeufm=eufm5
\newfam\eufmfam
\textfont\eufmfam=\teneufm
\scriptfont\eufmfam=\seveneufm
\scriptscriptfont\eufmfam=\fiveeufm
\def\mathfrak#1{{\fam\eufmfam\relax#1}}

\font\tenmsb=msbm10
\font\sevenmsb=msbm7
\font\fivemsb=msbm5
\newfam\msbfam
\textfont\msbfam=\tenmsb
\scriptfont\msbfam=\sevenmsb
\scriptscriptfont\msbfam=\fivemsb
\def\Bbb#1{{\fam\msbfam #1}}

\def \NN {\Bbb N}
\def \CC {\Bbb C}

\def \ZZ {\Bbb Z}

\def\rightheadline{{\hfil{\ff
The fourth moment of the zeta-function}\hfil\tenrm\folio}}

\def\leftheadline{{\tenrm\folio\hfil{\ff
A. Ivi\'c }\hfil}}
\def\emptyheadline{\hfil}
\headline{\ifnum\pageno=1 \emptyheadline\else
\ifodd\pageno \rightheadline \else \leftheadline\fi\fi}

\topmatter
\title THE FOURTH MOMENT OF THE ZETA-FUNCTION 
\endtitle
\author   Aleksandar Ivi\'c \endauthor
\address{ \bigskip
Aleksandar Ivi\'c, Katedra Matematike RGF-a
Universiteta u Beogradu, \DJ u\v sina 7, 11000 Beograd,
Serbia (Yugoslavia). \bigskip}
\endaddress
\email {\tt aleks\@ivic.matf.bg.ac.yu, 
aivic\@rgf.rgf.bg.ac.yu} \endemail
\abstract
{An overview of results and problems concerning the
asymptotic formula for $\int_0^T|\zt|^4\d t $ is given,
together with a discussion of modern methods from spectral
theory used in recent work on this subject.}
\endabstract
\endtopmatter

\noindent

\def\hf{{\textstyle{1\over2}}}
\def\txt#1{{\textstyle{#1}}}
\def\a{\alpha}
\def\e{\varepsilon}
\def\f{\varphi}
\def\G{\Gamma}
\def\H{H_j^3(\hf)}

\def\s{\sigma}

\def\={\;=\;}
\def\zx{\zeta(\hf+ix)}
\def\zt{\zeta(\hf+it)}

\font\teneufm=eufm10
\font\seveneufm=eufm7
\font\fiveeufm=eufm5
\newfam\eufmfam
\textfont\eufmfam=\teneufm
\scriptfont\eufmfam=\seveneufm
\scriptscriptfont\eufmfam=\fiveeufm
\def\mathfrak#1{{\fam\eufmfam\relax#1}}

\font\tenmsb=msbm10
\font\sevenmsb=msbm7
\font\fivemsb=msbm5
\newfam\msbfam
      \textfont\msbfam=\tenmsb
      \scriptfont\msbfam=\sevenmsb
      \scriptscriptfont\msbfam=\fivemsb
\def\Bbb#1{{\fam\msbfam #1}}

\def \NN {\Bbb N}
\def \CC {\Bbb C}

\def \ZZ {\Bbb Z}
\def  \HH {\Bbb H}
 
 \def\D{\Delta} \def\G{\Gamma} \def\e{\varepsilon}
\def\no{\noindent} 
\def\DJ{\leavevmode\setbox0=\hbox{D}\kern0pt\rlap
 {\kern.04em\raise.188\ht0\hbox{-}}D}

\heading 1. Introduction
\endheading
  A central place in analytic number theory is occupied by the Riemann
  zeta-function $\zeta(s)$, defined for $\R s > 1$ by
  $$
  \zeta(s) = \sum_{n=1}^\infty n^{-s} =
  \prod_{p \,\text {prime}} (1 - p^{-s})^{-1},
  $$
  and otherwise by analytic continuation. It admits meromorphic
  continuation to the whole complex plane, its only singularity
  being the simple pole $s = 1$ with residue 1. 
   From the functional equation (see [9], [32])
  $$
  \zeta(s) = \chi(s)\zeta(1 - s),\quad \chi(s) = 2^s\pi^{s-1}\sin\bigl(
  {\hf\pi s}\bigr)\Gamma(1 - s), 
  $$
  which is valid for any complex $s$, it follows that $\zeta(s)$
 has zeros at $s = -2, -4, \ldots\,$ . These zeros are traditionally
 called the ``trivial" zeros of $\zeta(s)$, to distinguish them
 from the complex zeros of $\zeta(s)$, of which the smallest ones
 (in absolute value) are ${1 \over 2} \pm 14.134725\ldots i$.
 It is well-known that all complex zeros of  $\zeta(s)$  lie
 in the so-called ``critical strip" $0 < \sigma = \R s < 1$,
  and the
Riemann Hypothesis (RH for short)  is the conjecture,
stated in 1859 by B. Riemann [31], that
{\it very likely all complex zeros of  $\zeta(s)$
have real parts equal to $\hf$}. For this reason the line
$\sigma = \hf$ is called the ``critical line" in the theory of  
$\zeta(s)$. The   RH is extensively discussed by many authors,
and recently  by E. Bombieri [3].

The RH is now probably the most celebrated and difficult open
problem in whole Mathematics. Its proof (or disproof) would
have very important consequences in multiplicative number theory,
especially in problems involving the distribution of primes.

\bigskip 
\heading
2. Mean values of $|\zt|$
\endheading
Mean values of $\vert\zeta({1\over 2} + it)\vert$
are fundamental in the theory and applications of $\z(s)$ (see
[11] for an extensive account). For $k \geqslant 1$ a fixed integer let
$$
\int_0^T\vert\zeta(\txt{1\over 2} + it)\vert^{2k}\d t =
T\,P_{k^2}(\log T) + E_k(T),\leqno(2.1)
$$
where for some suitable coefficients $a_{j,k}$ one has
$$
P_{k^2}(y) = \sum_{j=0}^{k^2}a_{j,k}y^j, \leqno(2.2) 
$$
and in particular
$$
P_1(y) = y + 2\gamma - 1 - \log(2\pi), 
$$
where $\gamma = -\G'(1) = 0.577\ldots\,$ is Euler's constant. One hopes that
$$
E_k(T) = o(T) \qquad (T \to \infty) 
$$
will hold for every fixed integer $k \geqslant 1$, but so far this is
known to be true only in the cases $k = 1$ and $k = 2$, when
$E_k(T)$ is a true error term in (2.1). When $k \geqslant 3$ it is not 
even known what should be the values of the coeeficients $a_{j,k}$ in (2.2).
The connection between $E_k(T)$ and the RH is indirect,
namely there is a connection with the {\it Lindel\"of hypothesis}
(LH for short). The LH is also a famous unsettled problem, and it
states that ($f \ll g$ means $|f(x)| < Cg(x)$ for some $C > 0$
and $x \geqslant x_0$)
$$
\zeta(\hf + it)  \ll_\e t^\varepsilon \leqno(2.3)
$$
for any given $\varepsilon > 0$ and $t \geqslant t_0 > 0$ (since
$\overline {\zeta({1\over 2} + it)} = \zeta({1\over 2} - it)$,
$t$ may be assumed to be positive). The
RH implies LH, in fact it even gives  (see [32]) an estimate stronger
than (2.3), namely
$$
\zeta(\hf + it) \ll 
\exp\left({A\,\log t \over \log \log t}
\right) \qquad (A > 0, t \geqslant t_0).
$$
It is yet unknown whether the LH  implies the RH.
The best current unconditional bound for the order of $\zt$, which is
far from the LH, 
is due to M.N. Huxley [7]. This is
$$
\zeta(\hf + it) \ll_\varepsilon t^{c+\varepsilon}, \qquad
c = {\txt{32\over205}} = 0.156098\ldots\,.
$$
Huxley  [7] also proved
$$
E_1(T) \,\ll_\e\, T^{{137\over432}+\e},\quad
{\txt{137\over432}} = 0.31713\ldots\,.
$$
\medskip \no
The LH is equivalent to the bound
$$
\int_0^T\vert\zeta(\hf + it)
\vert^{2k}\d t \ll_{k,\varepsilon}
T^{1+\varepsilon}  
$$
for every $k \geqslant 1$ and any $\varepsilon > 0$, which in 
turn is the same as
$$
E_k(T) \ll_{k,\varepsilon} T^{1+\varepsilon}\qquad(k \in \NN). 
$$
We have $\Omega$-results in the case $k = 1,2$, which show that 
$E_1(T)$ and $E_2(T)$ cannot be always small.  J.L. Hafner and 
the author [4] proved that
$$
E_1(T) = \Omega_+\bigl((T\,\log T)^{1\over 4}(\log \log T)^{3+\log4\over 4}
{\roman e}^{-C\sqrt{\log \log \log T}}\bigr) 
$$ \noindent
and
$$
E_1(T) = \Omega_-\biggl(T^{1\over 4}\exp\Bigl({D(\log \log T)^{1\over 4}
\over (\log \log \log T)^{3\over 4}}\Bigr)\biggr) 
$$
for some absolute constants $C,D >0$. Moreover, the author [10]
proved that there exist constants $A,B > 0$ such that, for $T \geqslant T_0$,
every interval $[T, T + B\sqrt{T}]$ contains points $t_1, t_2$ 
for which
$$
E_1(t_1) > At_1^{1\over4}, \quad E_1(t_2) < -At_2^{1\over4}.
$$
($f = \Omega_+(g)$ means that  $\limsup f/g > 0$,
$f = \Omega_-(g)$ means that $\liminf f/g < 0$).

\bigskip 
\heading
3. The fourth moment of $|\zt|$
\endheading
The asymptotic formula for the fourth moment of $\z(s)$ on the 
critical line is customarily written as (this is (2.1) and (2.2) when $k=2$)
$$
\int_0^T|\zt|^4\d t \;=\; TP_4(\log T) \;+\;E_2(T), \quad
P_4(x) \;=\; \sum_{j=0}^4\,a_jx^j.
$$
A classical result of A.E. Ingham  [8] from 1926
is that $a_4 = 1/(2\pi^2)$ and that the error term $E_2(T)$
satisfies the bound $E_2(T) \ll T\log^3T$. In 1979 D.R. Heath-Brown [6]
made significant progress in this problem by proving that 
$$
E_2(T) \ll_\e T^{{7\over8}+\e}.
$$
He also calculated
$$
a_3 \;=\; 2(4\gamma - 1 - \log(2\pi) - 12\zeta'(2)\pi^{-2})\pi^{-2}
$$
and produced more complicated expressions for $a_0, a_1$ and $a_2$.
The author [12] made an explicit evaluation of
$a_0, a_1,a_2$. For $\R s > 1$
$$
\zeta^2(s) = \sum_{n=1}^\infty d(n)n^{-s},
$$
where $d(n)$ is the sum of positive divisors of $n$. Note that
$\z^4 = \z^2\cdot\z^2$ and that $\z^2$ can be approximated by a finite
sum $\sum_n d(n)n^{-s}$. Therefore by integrating 
$$
|\zt|^4 = \zeta^2(\hf+it)\zeta^2(\hf-it)
$$
we are led to the asymptotic evaluation of the sum
$$
\sum_{n\leqslant x}d(n)d(n + f),
$$
where $1 \leqslant f \leqslant x$ is { not} fixed. This is the so-called
binary additive divisor problem (see [20], [22] and [27]). Modern
approaches both to the study of $E_2(T)$ and the
binary additive divisor problem  involve the use of spectral theory.

\bigskip 
\heading
4. Spectral theory
\endheading
\medskip\no
For a competent account on spectral theory and its applications to
$\z(s)$ the reader is referred to Y. Motohashi's monograph [30]. Here
we shall only briefly present some basic facts.
On the upper complex half-plane $\HH$ the modular group
$$
\G = {\roman SL}(2,\ZZ) = \left\{{a\quad b\choose c\quad d}  \;:\;
(a,b,c,d\in \ZZ) \wedge (ad -b c = 1)\right\} 
$$ 
acts by $\gamma z = (az + b)/(cz + d)$ if $\gamma \in \G$.
The { non-Euclidean Laplace operator}
$$
L = -y^2\left({\partial^2\over\partial x^2} +
{\partial^2\over\partial y^2} \right)
$$
is invariant under $\G$, as is the measure
$
\d\mu(z) = y^{-2}\d x\d y\, (z = x + iy).
$
Non-holomorphic cusp forms (so-called Maass wave forms) are 
eigenfunctions  $\Psi(z)$
of the discrete spectrum of $L$, which has the form
$$
\left\{\lambda_j\right\}_{j=1}^\infty,\qquad \lambda_j 
= \kappa^2_j + {\txt{1\over4}}\quad(\kappa_j > 0).
$$
They satisfy $L\Psi(z) = \lambda\Psi(z)$, $\Psi(\gamma z) = \Psi(z)$ 
for $\gamma \in \G$, and the finiteness condition
$$
\int_{\Cal D} |\Psi(z)|^2\d\mu(z) = \int_{\Cal D}
|\Psi(x+iy)|^2y^{-2}\d x\d y \;< \; \infty.
$$
Here ${\Cal D}$ is the fundamental domain of $\G$, namely
$$
{\Cal D}  = \left\{z\;:\; y > 0, |z| > 1,\,
-\hf \leqslant x \leqslant \hf
\right\} \cup  \left\{z\;:\; |z| = 1,  -\hf \leqslant x \leqslant 0
\right\}. 
$$
Let henceforth $\f_j$ be the Maass wave form attached to $\kappa_j$
so that $\{\f_j\}_{j=1}^\infty$ forms an orthonormal basis with
respect to the { Petersson inner product}
$$
(f_1,\,f_2) := \int_{\Cal D}f_1{\bar f_2}\d\mu(z),
$$
and $\f_j$ is an eigenfunction of every { Hecke operator}.
The Hecke operator $T_n$ acts on $\HH$, for given $n\in\NN$, by
the relation
$$
(T_nf)(z) = n^{-1/2}\sum_{ad=n,d>0}\;\sum_{b(\text{mod}\,d)}
f\left({az+b\over d}\right).
$$
We have the Fourier expansion
$$
\f_j(z) = \sum_{n\not=0}\rho_j(n)e^{2\pi ix}\sqrt{y}K_{i\kappa_j}
(2\pi|n|y)$$
for $z = x+iy,\,\rho_j(n) = \overline{\rho_j(-n)},$
where $K$ is the { Bessel function (the Mcdonald function)}
$$
K_s(z) = \hf\int_0^\infty t^{s-1}\exp\left(-\hf z(t + {1\over t})
\right)\d t\quad(\R z > 0).
$$
For $n\in\NN$, $t_j(n)$ is the eigenvalue corresponding to $\f_j$
with respect to $T_n$,
$$
T_n\f_j(z) = t_j(n)\f_j(z).
$$
Then
$$
\rho_j(1)t_j(n) \;=\;\rho_j(n).
$$
The  Hecke series attached to $\f_j(z)$ is
$$
H_j(s) := \sum_{n=1}^\infty t_j(n)n^{-s} 
= \prod_p(1 - t_j(p)p^{-s} + p^{-2s})^{-1}\quad(\R s > 1).
$$
It is known that $H_j(s)$ continues analytically to an entire function 
over $\CC$, and that for any $s$ satisfies the functional equation
$$
H_j(s) = \pi^{-1}(2\pi)^{2s-1}\G(1-s+i\kappa_j)\G(1-s-i\kappa_j)
\Bigl\{-\cos(\pi s) + \e_j\cosh(\pi\kappa_j)\Bigr\}H_j(1-s)
$$
with $\e_j = \pm 1$. 

\bigskip 
\heading
5. Motohashi's explicit formula
\endheading
\medskip\no
Recent progress on the fourth moment of $|\zt|$ is primarily
due to the fundamental explicit formula of Y. Motohashi
(see [11], [25], [26], [30]) of 1989 for
$$ 
I(T,\D) := {1\over\D\sqrt{\pi}}\int_{-\infty}^\infty
|\zeta(\hf+iT+it)|^4{\roman e}^{-t^2\D^{-2}}\d t,
$$
under the condition
$$
0 < \D \leqslant {T\over\log T}.
$$
The presence of the  Gaussian smoothing factor ${\roman e}
^{-x^2}$ enabled Motohashi to deal with various convergence
problems occurring in the proof. In simplified form,
the formula is
$$\eqalign{
I(T,\D) &= \pi 2^{-1/2}T^{-1/2}\sum_{j=1}^\infty
\a_j \kappa_j^{-1/2}H^3_j(\hf)\sin\left(\kappa_j\log\left(
{\kappa_j\over4{\roman e} T}\right)\right){\roman e}^
{-(\D\kappa_j/2T)^2}\cr& + O(\log^CT).  \cr}\leqno(5.1)
$$
Here
$$
\a_j \;=\;{|\rho_j(1)|^2\over\cosh(\pi\kappa_j)},\leqno(5.2)
$$
and,  for arbitrary fixed $A > 0$, $C = C(A) > 0$ with $\D$ satisfying
$$
{\sqrt{T} \over\log^AT} \leqslant \D \leqslant T\exp(-\sqrt{\log T}).
$$
The proof of the formula for $I(T,\D)$ depends heavily on the
use of spectral theory, especially the so-called Kuznetsov
trace formulas, which relate sums of Kloosterman sums (see (6.4)) to certain
sums involving quantities from spectral theory.
\bigskip 
\heading
6. New results on the fourth moment
\endheading

\no At present the bound
$$
\int_0^T\,|\zt|^k\d t \;\ll_\e\; T^{1+\e}
$$
is not known to hold  for any constant $k > 4$. Therefore
the function $E_2(T)$ is particularly important in the theory
of mean values of $\z(s)$. In recent years
much advance has been made in connection with
$E_2(T)$ and related problems. One of the main problems is to get rid of
the  Gaussian smoothing factor ${\roman e}
^{-x^2}$ in Motohashi's formula (5.1) and apply it to obtain
results on $E_2(T)$ itself.  

\medskip\no 
Motohashi and the author in four papers [17]--[20] obtained
several results on $E_2(T)$ and some related problems. It is known now that
($f = \Omega (g)$ means that $\limsup|f|/g > 0$)
$$
E_2(T) \;=\; O(T^{2/3}\log^{C_1}T),\quad
 E_2(T) \;=\; \Omega(T^{1/2}),
$$
$$
\int\limits_0^TE_2(t)\d t \;=\; O(T^{3/2}),
\;\int\limits_0^TE_2^2(t)\d t
\;=\; O(T^2\log^{C_2}T),$$
with effective constants $C_1,\,C_2 > 0$ 
(the values $C_1 = 8, C_2 = 22$ are
admissible). 
Y. Motohashi [28] improved the omega-result to 
$$
E_2(T) = \Omega_\pm(T^{1/2}).
$$
Finally  the author [15] made further progress in this problem
by proving the following
quantitative omega-result: there exist two 
constants $A >0,\,B > 1$ such that
for $T \geqslant T_0 > 0$ every interval $[T,\,BT]$
contains points $T_1,T_2$ for which
$$
E_2(T_1) \;>\; AT_1^{1/2},\quad E_2(T_2) \;<\; -AT_2^{1/2}.
$$
This follows from the asymptotic formula
$$\eqalign{
\int\limits_0^\infty E_2(t){\roman e}^{-t/T}\d t &=
2T^{3\over2}\,\R \left\{\sum_{j=1}^\infty \a_jH_j^3(\hf)R(\kappa_j)\G(\hf
 {-i\kappa_j})T^{-\kappa_j} \right\} \cr&
+ O\left(T^{{3\over2}}\exp\left(-{\log T\over
\log\log T}\right)\right), \cr}\leqno(6.1)
$$
with  $\a_j$ given by (5.2) and
$$
R(y) \;:=\;\sqrt{{\pi\over2}}{\Bigl(2^{iy}{\G({1\over4}
+ {i\over2}y)
\over\G({1\over4} - {i\over2}y)}\Bigr)}^3\G(-2iy)
{\cosh}(\pi y).\leqno(6.2)
$$
Y. Motohashi and the author [19] proved that
$$
\sum_{r=1}^R\int_{t_r}^{t_r+\D}|\zt|^4\d t \;\ll\;
R\D\log^4T + TR^{1/2}\D^{-1/2}\log^CT
$$
for some $C > 0$, where $\log T \;\ll\; \D \ll {T\over\log T}$ and
$$
T \leqslant t_1 < t_2 < \ldots < t_R \leqslant 2T,\quad
t_{r+1} - t_r \geqslant \D \quad(r = 1,\ldots,R-1).
$$
From this result it is not difficult to obtain as a corollary the bound
$$
\int_0^T|\zt|^{12}\d t \;=\; O(T^2\log^BT),
$$
proved first (with $B = 17$) by Heath-Brown [5] in 1978.

\smallskip
The asymptotic formula (6.1) for
$\int_0^\infty E_2(t){\roman e}^{-t/T}\d t $ is a
Laplace transform formula. Integral transforms such as the Mellin
transform (see [13], [21]) and the Laplace transform (see [16], [23])
play an important r\^ole in analytic number theory, in particular in the
theory of $\z(s)$. One can consider the general function
$$
L_k(s) \;:=\; \int\limits_0^\infty |\zx|^{2k}e^{-sx}\d x
\qquad(k \in \NN,\, \R s > 0).
$$
A classical result of H. Kober [24] from 1936
says that, as $\s \to 0+$,
$$
L_1(2\s) = {\gamma-\log(4\pi\s)\over2\sin\s} +
\sum_{n=0}^Nc_n\s^n + O(\s^{N+1})
$$
for any given integer $N \geqslant 1$, where the $c_n$'s are effectively
computable constants.

F.V. Atkinson [2] obtained in 1941 the asymptotic formula
$$
L_2(\s) =  
{1\over\s}\left(A\log^4{1\over\s} + B\log^3{1\over\s} + C\log^2{1\over\s} 
+ D\log {1\over\s} + E\right) + \lambda_2(\s),\leqno(6.3)
$$
where $\s \to 0+$,
$$
A = {1\over2\pi^2},\,B =\pi^{-2}(2\log(2\pi) -
 6\gamma + 24\z'(2)\pi^{-2}),\;
 \lambda_2(\s) \;\ll_\e\;\left({1\over\s}\right)^{{13\over14}+\e}.
$$
He also indicated how, by the use of estimates for  Kloosterman sums
$$
S(m,n;c) := \sum_{1\leqslant d < c,(d,c)=1,dd'\equiv
 1(\text{mod}\,c)}{\roman e}\left({md+nd'\over c}\right)
\quad({\roman e}(z) = {\roman e}^{2\pi i z}), \leqno(6.4)
$$
one can improve the exponent ${13\over14}$  to ${8\over9}$. This is
important historically, in view of  contemporary importance of
 Kloosterman sums which stems from the Kuznetsov trace formulas and
other important applications.

The author [12] gave explicit, albeit complicated expressions for
the remaining coefficients $C,D$ and $E$ in (6.3) and improved
Atkinson's bound for    $\lambda_2(\s)$ to
$$
\lambda_2(\s) \;\ll\; \s^{-1/2}\qquad(\s\to 0+).
$$
Recently in [16] he proved a result which generalizes and sharpens (6.3):
Let $0 \leqslant \phi < {\pi\over2}$ be given. Then for $0 < |s| \leqslant 1$ and
$|\arg s| \leqslant \phi$ we have
$$\eqalign{
L_2(s)  &
= {1\over s}(A\log^4{1\over s} + B\log^3{1\over s} \;+
C\log^2{1\over s} + D\log{1\over s} + E) \cr&
+ \,s^{-{1\over2}}\Bigl\{\sum_{j=1}^\infty \a_j\H\Bigl(
s^{-i\kappa_j}R(\kappa_j)\G(\kappa_j)
+ s^{i\kappa_j}R(-\kappa_j)\G(-\kappa_j) \Bigr)\Bigr\}
+ G_2(s),\cr}
$$
where the constants $A,\ldots\,,E$ are as in (6.3), $R(y)$ is given 
by (6.2) and in the above region $G_2(s)$ is a regular function satisfying
($C > 0$ is a suitable constant)
$$ 
G_2(s)  
\ll |s|^{-{1\over2}}\exp\left\{
-{C\log(|s|^{-1})\over(\log\log(|s|^{-1}))^{2\over3}
(\log\log\log(|s|^{-1}))^{1\over3}}\right\}.
$$
In [14] the author proved:
There exist constants $A > 0$ and $B > 1$ such that,
for $T \geqslant T_0 > 0$, every interval $\,[T,\,BT]\,$ contains points
$t_1,\,t_2$ for which
$$
\int_0^{t_1}E_2(t)\d t \;>\; At_1^{3/2},\quad
\int_0^{t_2}E_2(t)\d t \;<\; -At_2^{3/2}.
$$

\medskip\no
This result implies that 
$$
\int_0^T E_2(t)\d t \;=\;\Omega_\pm(T^{3/2}),
$$
and that
$$
\limsup_{n\to\infty}\,{\log(u_{n+1}-u_n)\over\log u_n} \;\leqslant \;1,
$$
where $u_n$ is the $n$--th zero of $E_2(T)$. It can also be used 
to prove the lower bound result
$$
\int_0^TE_2^2(t)\d t \;\gg\;T^2,
$$
which complements the earlier result with Motohashi (cf. [18]) that
$$
\int_0^TE_2^2(t)\d t \;\ll\;T^2\log^CT\qquad(C = 22). 
$$
Another recent result of the author [15], which provides 
an asymptotic formula for the integral of $E_2(t)$, is
$$\eqalign{&
\int_0^T E_2(t) \d t =  
O\{T^{3\over2}\exp(-C(\log T)^{3\over5}(\log\log T)^{-{1\over5}})\} \;+\cr&
+\, 2T^{3\over2}\R\left\{
\sum_{j=1}^\infty\a_j\H{T^{i\kappa_j}\over({\txt{1\over2}}+i\kappa_j)
({\txt{3\over2}}+i\kappa_j)}R(\kappa_j)\right\}.\cr}
$$
Here $C>0$ is a suitable constant, and the error term depends on the
best known zero-free region for $\z(s)$ (see [9, Chapter 6]).

\medskip
In spite of significant recent results on $E_2(T)$, many open
problems remain. Here are four of them.

\medskip\no
{\bf Problem 1.} Does there exist $A > 0$ such that, as $T\to\infty$, 
$$
\int_0^TE_2^2(t)\d t \;\sim\;AT^2 ?
$$

\medskip\no
{\bf Problem 2.} Is it true that for every $\e > 0$
$$
E_2(T) \;=\; O_\e\left(T^{{1\over2}+\e}\right)?
$$

\medskip\no
{\bf Problem 3.} Does one have
$$
\limsup_{T\to\infty} |E_2(T)|T^{-{1\over2}} \;=\; \infty?
$$
\medskip\no
{\bf Problem 4.} Does one have
$$
\limsup_{n\to\infty}\,{\log(u_{n+1}-u_n)\over\log u_n} \;<\;1?
$$

\vfill
\eject
\topglue2cm

\Refs

\item{[1]} F.V. Atkinson, `The mean value of the zeta-function on
the critical line', {\it Quart. J. Math. Oxford} {\bf 10}(1939), 122-128.

\item{[2]} { F.V. Atkinson}, `The mean value of the zeta-function on the
critical line', {\it Proc. London Math. Soc.} (2){\bf 47}(1941), 174-200.

\item{[3]} E. Bombieri, `Problems of the Millenium: the Riemann Hypothesis',
2000, {\tt http//www.ams.org\ /claymath/prize{\_}problems/riemann.pdf.},
11pp.

\item{ [4]} J.L. Hafner and A. Ivi\'c, `On the mean square of the
Riemann zeta-function on the critical line', {\it J. Number Theory}
{\bf 32}(1989), 151-191.

\item{[5]}  D.R. Heath-Brown, `The twelfth power moment of the
Riemann zeta-function', {\it Quart. J. Math. Oxford} {\bf29}(1978), 443-462.

\item{[6]} { D.R. Heath-Brown}, `The fourth moment of the Riemann 
zeta-function',
{\it Proc. London Math. Soc.}  (3){\bf38}(1979), 385-422.

\item{[7]} M.N. Huxley, `Integer points, exponential sums and the
Riemann zeta-function', Proc. of ``Millenial Conference on Number
Theory", Urbana, May 2000, in print.

\item{[8]} { A.E. Ingham}, `Mean-value theorems in the theory of 
the Riemann
zeta-function', {\it Proc. London Math. Soc.}  (2){\bf27}(1926), 273-300.

\item {[9]} { A. Ivi\'c},  `The Riemann zeta-function', {\it John Wiley
and Sons}, New York, 1985.

\item{[10]} {A. Ivi\'c}, `Large values of certain 
number-theoretic error terms', {\it Acta Arith.} {\bf56} (1990), 135-159.

\item {[11]} { A. Ivi\'c},  `Mean values of the Riemann zeta-function',
LN's {\bf82}, {\it Tata Institute of Fundamental Research}, Bombay, 1991 
(distr. by Springer Verlag, Berlin etc.).

\item{ [12]} { A. Ivi\'c},  `On the fourth moment of the Riemann
zeta-function', {\it Publs. Inst. Math. (Belgrade)} 
{\bf57(71)}(1995), 101-110.

\item{[13]} { A. Ivi\'c},  `The Mellin transform and 
the Riemann zeta-function', in
{\it ``Proceedings of the Conference on Elementary and
Analytic Number Theory 
(Vienna, July 18-20, 1996)"}, Universit\"at Wien \& Universit\"at f\"ur
Bodenkultur, Eds. W.G. Nowak and J. Schoi{\ss}engeier, Vienna 1996, 112-127.

\item{[14]} A. Ivi\'c, `On the error term for the fourth moment of the
Riemann zeta-function', {\it J. London Math. Soc.}, 
{\bf60}(2)(1999), 21-32.

\item{[15]} A. Ivi\'c,  `On the integral of the error term in 
the fourth moment of the Riemann zeta-function', {\it Functiones 
et Approximatio}, {\bf28}(2000), 37-48.

\item{[16]} A. Ivi\'c, `The Laplace transform of the fourth moment
of the zeta-function', {\it Univ. Beog. Publik. Elektroteh. Fak. Ser.
Matematika} {\bf11}(2000), 41-48.

\item{ [17]}  A. Ivi\'c and Y. Motohashi,  `A note on the mean value of
the zeta and L-functions VII', {\it  Proc. Japan Acad. Ser. A}
{\bf  66}(1990), 150-152.

\item{ [18]}{ A. Ivi\'c and Y. Motohashi}, `The mean square of the
error term for the fourth moment of the zeta-function', 
{\it Proc. London Math.
Soc.} (3){\bf66}(1994), 309-329.

\item {[19]} { A. Ivi\'c and Y. Motohashi},  `The fourth moment of the
Riemann zeta-function', {\it J. Number Theory} {\bf 51}(1995), 16-45.

\item {[20]} { A. Ivi\'c and Y. Motohashi},  `On some estimates
involving the binary additive problem', {\it Quart. J. Math. Oxford}
(2){\bf46}(1995), 471-483.

\item {[21]}  A. Ivi\'c, M. Jutila and Y. Motohashi,
{\it The Mellin transform of powers of the zeta-function}, {\it Acta
Arith.} {\bf95}(2000), 305-342.

\item{ [22]} { M. Jutila}, `The fourth moment of Riemann's zeta-function 
and the
additive divisor problem', in {\it ``Analytic Number Theory; Proc. Conf. in
Honor of Heini Halberstam, Vol. 2"} (eds. B.C. Berndt et al.), Birkh\"auser, 
Boston etc., 1996, 517-536.

\item{[23]} { M. Jutila}, `Mean values of Dirichlet series via Laplace 
transforms',
in {\it ``Analytic Number Theory"} (ed. Y. Motohashi), London Math. Soc.
LNS 247, {\it Cambridge University Press}, Cambridge, 1997, 169-207.

\item {[24]} H. Kober, Eine Mittelwertformel der Riemannschen Zetafunktion,
{\it Compositio Math.} {\bf3}(1936), 174-189.

\item {[25]} Y. Motohashi,  The fourth power mean of the Riemann
zeta-function, in {\it Proceedings of  the Amalfi Conference on
Analytic Number Theory 1989}, eds.  E. Bombieri et al., Universit\`a
di Salerno, Salerno, 1992, 325-344.

\item{ [26]} Y. Motohashi,   `An explicit formula for the fourth power
mean of the Riemann zeta-function', {\it Acta Math. }{\bf 170}(1993),
181-220.

\item {[27]} Y. Motohashi,`The binary additive divisor problem', {\it
Ann. Sci. \'Ecole. Norm. Sup. $4^e$ s\'erie} {\bf27}(1994), 529-572.

\item {[28]} Y. Motohashi,  `A relation  between the Riemann
zeta-function and the hyperbolic Laplacian', {\it Annali Scuola Norm.
Sup. Pisa, Cl. Sci. IV ser.} {\bf 22}(1995), 299-313.

\item {[29]} Y. Motohashi,  `The Riemann zeta-function and the
non-Euclidean Laplacian', {\it Sugaku Expositions}, AMS {\bf 8}(1995),
59-87.

\item {[30]} Y. Motohashi,  `Spectral theory of the Riemann
zeta-function', {\it Cambridge University Press}, Cambridge, 1997.

\item {[31]} B. Riemann, `\"Uber die Anzahl der Primzahlen
unter einer gegebener Gr\"osse', {\it Monatsber. Akad. Berlin}
(1859), 671-680.

\item {[32]} E.C. Titchmarsh,  The Theory of the Riemann Zeta-Function,
{\it Clarendon Press}, Oxford, 1951.

\bigskip
\bigskip
\bigskip

Aleksandar Ivi\'c

Katedra Matematike RGF-a

Universitet u Beogradu, \DJ u\v sina 7

11000 Beograd, Serbia (Yugoslavia)

\tt aivic\@rgf.bg.ac.yu, \enskip aivic\@matf.bg.ac.yu

\endRefs

\bye